\documentclass[10pt]{amsart}

\usepackage{amsmath}
\usepackage{amstext}
\usepackage{amssymb}
\usepackage{latexsym}
\usepackage{exscale}
\usepackage[english]{babel}
\usepackage{indentfirst}

\newcommand{\demo}{\vskip0.1cm \noindent {\it Proof. }}
\newcommand{\demodel}[1]{\vskip0.1cm \noindent {\it Proof of #1. }}
\renewcommand{\qed}{\hfill$\Box$}

\newcommand{\io}{\int_{\mathcal{{C}}_\O}}

\newcommand{\weig}{{y^{1-\a}}}

\renewcommand{\a}{\al}

\renewcommand{\d }{\delta }

\renewcommand{\O}{\Omega}

\newcommand{\dist}{\mbox{dist }}

\newcommand{\p}{\partial}

\newcommand{\be}{\beta}
\newcommand{\lan}{\lambda}
\newcommand{\Lan}{\Lambda}
\renewcommand{\l}{\lan}

\newcommand{\epsi}{\varepsilon}

\newcommand{\vfi}{\varphi}

\newcommand{\al}{\alpha}
\newcommand{\Om}{\Omega}
\newcommand{\om}{\omega}
\newcommand{\de}{\delta}
\renewcommand{\div}{\mbox{\rm div}}

\newcommand{\tr}{\mbox{\rm Traza}}

\newcommand{\grad}{\nabla}
\newcommand{\n}{\grad}

\newcommand{\dyle}{\displaystyle}

\newtheorem{thm}{Theorem}[section]
\newtheorem{lemma}[thm]{Lemma}
\newtheorem{prop}[thm]{Proposition}
\newtheorem{corol}[thm]{Corollary}

\newtheorem{defi}[thm]{Definition}
\newtheorem{obs}{Remark}[section]
\newenvironment{obse}{\begin{obs}\rm}{\end{obs}}

\def\diver{\mathop{\rm div}\nolimits}
\def\tr{\mathop{\rm Tr}\nolimits}
\def\ext{\mathop{\rm E}\nolimits}

\title[A concave-convex problem involving the fractional
Laplacian]{\bf A concave-convex elliptic problem involving the fractional
Laplacian}

\author[Br\"andle]{C. Br\"andle}

\author[Colorado]{E. Colorado}

\author[de Pablo]{A. de Pablo}
 \email{cbrandle@math.uc3m.es, ecolorad@math.uc3m.es, arturop@math.uc3m.es}
 \address{Departamento de Matem{\'a}ticas, U.~Carlos III de Madrid,
28911 Legan{\'e}s (Madrid), Spain}

\begin{document}
\begin{abstract}
We study a nonlinear elliptic problem defined in a bounded domain involving
fractional powers of the Laplacian operator together with a concave-convex
term. We characterize completely the range of parameters for which solutions of
the problem exist and prove a multiplicity result.
\end{abstract}

\maketitle

\section{Introduction}\label{sect-intro}

In the past decades the problem
$$
\left\{
    \begin{array}
      {rcl@{\qquad}l}
        -\Delta u&=&f(u)&\mbox{in }\Omega\subset\mathbb{R}^N,\\
    u&=&0&\mbox{on }\partial \Omega,
    \end{array}
    \right.
$$
has been widely investigated. See~\cite{A} for a survey, and for
example the list (far from complete)~\cite{ABC,Brezis-Nirenberg2,Lions85b} for more specific
problems, where different nonlinearities and different classes of domains,
bounded or not, are considered. Other different diffusion operators, like the $p$--Laplacian, fully nonlinear operators, etc, have been also treated,  see for
example~\cite{BEP, Caffarelli-Cabre,GP} and the references there in. We deal
here with a nonlocal version of the above problem, for a particular type of
nonlinearities, i.e., we study a concave-convex problem involving the fractional
Laplacian operator
\begin{equation}
  \label{prob:main.nonlocal}
\left\{
    \begin{array}
      {rcl@{\qquad}l}
        (-\Delta)^{\alpha/2}u&=&\lambda u^q+u^p,\qquad u>0&\mbox{ in } \Omega,\\
    u&=&0&\mbox{ on } \partial \Omega,
    \end{array}
    \right.
\end{equation}
with $0<\alpha<2$, $0<q<1<p< \frac{N+\alpha}{N-\alpha}$, $N>\alpha$,
$\lambda>0$ and $\Omega\subset\mathbb{R}^N$ a smooth bounded domain.

The nonlocal operator $(-\Delta)^{\alpha/2}$ in $\mathbb{R}^N$ is
defined on the Schwartz class through the Fourier transform,
$$
    [(-\Delta)^{\alpha/2}  g]^{\wedge}\,(\xi)=(2\pi|\xi|)^\alpha \,\widehat{g}(\xi),
$$
or via the Riesz potential, see for example~\cite{Landkof,Stein} for the
precise formula. As usual, $\widehat g$ denotes the Fourier Transform of $g$,
$\widehat g(\xi)=\int_{\mathbb{R}^N}e^{-2\pi ix\cdot\xi}g(x)\,dx$. Observe that
$\alpha=2$ corresponds to the standard local Laplacian.

This type of diffusion operators arises in several areas such as physics,
probability and finance, see for
instance~\cite{Applebaum,Bertoin,Cont-Tankov,VIKH}. In particular, the
fractional Laplacian can be understood as the infinitesimal generator of a stable L{\'e}vy
process,~\cite{Bertoin}.

There is another way of defining this operator. In fact, in the case
$\alpha=1$ there is an explicit form of calculating the half-Laplacian acting on a function $u$ in the
whole space $\mathbb{R}^N$, as the normal derivative on the boundary of its
harmonic extension to the upper half-space $\mathbb{R}^{N+1}_+$, the so-called
Dirichlet to Neumann operator. The ``$\alpha$ derivative"
$(-\Delta)^{\alpha/2}$ can be characterized in a similar way, defining the
$\alpha$-harmonic extension to the upper half-space,
see~\cite{Caffarelli-Silvestre} and Section~\ref{sect:whole-space} for details.
This extension is commonly used in the recent literature since it allows to
write nonlocal problems in a local way and this permits to use the variational techniques for these
kind of problems.

In the case of the operator defined in bounded domains $\Omega$, the above
characterization has to be adapted. The fractional powers of a linear positive
operator in $\Omega$ are defined by means of the spectral decomposition.
In~\cite{Cabre-Tan},  the authors consider the fractional operator
$(-\Delta)^{1/2}$ defined using the mentioned Dirichlet to Neumann operator,
but restricted to the cylinder
$\Omega\times\mathbb{R}_+\subset\mathbb{R}^{N+1}_+$, and show that this
definition is coherent with the spectral one, see also~\cite{stinga-torrea} for
the case $\alpha\neq1$. We recall that this is not the unique possibility of
defining a nonlocal operator related to the fractional Laplacian in a bounded
domain. See for instance the definition of the so called {\it regional
fractional Laplacian} in~\cite{BBC,Guan-Ma}, where the authors consider the
Riesz integral
 restricted to the domain $\Omega$. This leads to a different operator related to a Neumann problem.

As to the concave-convex nonlinearity, there is a huge amount of results
involving different (local) operators, see for
instance~\cite{AbCP,ABC,BEP,CCP,CP,GP}. We quoted the work~\cite{ABC} from where
some ideas are used in the present paper. In most of the problems considered in
those papers a critical exponent appears, which generically separates the range
where compactness results can be applied or can not (in the fully nonlinear
case the situation is slightly different, but still a critical exponent
appears,
\cite{CCP}). In our case, the critical exponent with respect to the corresponding Sobolev
embedding is given by $2^*_\alpha=\frac{2N}{N-\alpha}$. This is a reason why
problem \eqref{prob:main.nonlocal} is studied in the subcritical case
$p<2^*_\alpha-1=\frac{N+\alpha}{N-\alpha}$; see also the nonexistence result
for supercritical nonlinearities in Corollary~\ref{supercritical}.

The main results we prove characterize the existence of solutions
of~\eqref{prob:main.nonlocal} in terms of the parameter $\lambda$. A
competition between the sublinear and superlinear terms plays a role, which
leads to different results concerning existence and multiplicity of solutions, among others.
\begin{thm}\label{thm:main}
  There exists $\Lambda>0$ such that for Problem~\eqref{prob:main.nonlocal}
  there holds:
\begin{itemize}
\item[{\rm 1.}] If $0<\lambda <\Lambda$ there is a minimal solution.
Moreover, the family of minimal solutions is increasing with respect to
$\lambda$.
 \item[{\rm 2.}] If $\lambda =\Lambda$ there is at least one solution.
 \item[{\rm 3.}] If $\lambda >\Lambda$ there is no solution.
\end{itemize}
\end{thm}
Moreover, we have the following multiplicity result:
\begin{thm}\label{thm:main2}
For each $0<\lambda <\Lambda$, Problem~\eqref{prob:main.nonlocal} has at least two solutions.
\end{thm}
On the contrary, there is a unique solution with small norm.
\begin{thm}
\label{thm:smallnorm}
For each $0<\lambda <\Lambda$ fixed,  there exists a constant $A>0$ such that there exists
  at most one solution $u$
  to Problem~\eqref{prob:main.nonlocal} verifying $$\|u\|_\infty\le A.$$
\end{thm}

For $\alpha\in[1,2)$ and $p$ subcritical, we also prove that there exists an
universal $L^\infty$-bound for every solution independently of $\lambda$.
\begin{thm}\label{th:universalbound}
  Let $\alpha\ge1$. Then there exists a constant $C>0$ such that, for any
  $0<\lambda\le\Lambda$,
  every solution to Problem~\eqref{prob:main.nonlocal} satisfies
  $$
\|u\|_{\infty}\le C.
  $$
\end{thm}
The proof of this last result relies on the classical argument of rescaling
introduced in \cite{Gidas-Spruck} which yields to problems on unbounded domains, which require some Liouville-type results, which can be seen in~\cite{dps}. This is the point where the
restriction $\alpha\ge1$ appears.

The paper is organized as follows: in Section~\ref{sect:whole-space} we
recollect some  properties of the fractional Laplacian in the whole space, and
establish a trace inequality corresponding to this operator; the fractional
Laplacian in a bounded domain is considered in Section~\ref{sect:bounded}, by
means of the use of the $\alpha$-harmonic extension; this includes studying an
associated linear equation in the local version. The main section,
Section~\ref{sect-f(u)}, contains the results related to the nonlocal nonlinear
problem~\eqref{prob:main.nonlocal}, where we prove
Theorems~\ref{thm:main}--\ref{th:universalbound}.


\section{The fractional Laplacian in $\mathbb{R}^{N}$}
\label{sect:whole-space}
\setcounter{equation}{0}

\subsection{Preliminaries}

Let $u$ be a regular function in $\mathbb{R}^N$. We say that
$w=\ext_\alpha(u)$ is the $\alpha$-harmonic extension of $u$ to the
upper half-space, $\mathbb{R}^{N+1}_+$, if $w$ is a solution to the
problem
$$
  \left\{
    \begin{array}
      {rcl@{\qquad}l}
        -\diver(y^{1-\alpha}\nabla w)&=&0&\mbox{ in } \mathbb{R}^{N+1}_+,\\ [5pt]
        w&=&u&\mbox{ on } \mathbb{R}^N\times\{y=0\}.
    \end{array}
    \right.
$$
In~\cite{Caffarelli-Silvestre} it is proved that
 \begin{equation}
\lim\limits_{y\to0^+}y^{1-\alpha}\dfrac{\partial w}{\partial
y}(x,y)=-\kappa_\alpha(-\Delta)^{\alpha/2}u(x),
\label{normalder}\end{equation} where
$\kappa_\alpha=\frac{2^{1-\alpha}\Gamma(1-\alpha/2)}{\Gamma(\alpha/2)}$.
Observe that $\kappa_\alpha=1$ for $\alpha=1$ and
$\kappa_\alpha\sim1/(2-\alpha)$ as $\alpha\to2^-$. As we pointed out
in the Introduction, identity~\eqref{normalder} allows to formulate
nonlocal problems involving the fractional powers of the Laplacian
in $\mathbb{R}^N$ as local problems in divergence form in the
half-space $\mathbb{R}^{N+1}_+$.

The appropriate functional spaces to work with are
$X^\alpha(\mathbb{R}^{N+1}_+)$ and $\dot{H}^{\alpha/2}(\mathbb{R}^N)$, defined
as the completion of ${C}_0^\infty(\overline{\mathbb{R}^{N+1}_+})$ and
$C_0^\infty(\mathbb{R}^N)$, respectively, under the norms
$$
 \begin{array}{l}
    \displaystyle\|\phi\|^2_{X^\alpha}=\int_{\mathbb{R}^{N+1}_+} y^{1-\alpha}|\nabla
    \phi(x,y)|^2\,dx dy,\\[12pt]
    \displaystyle\|\psi\|^2_{\dot{H}^{\alpha/2}}=\int_{\mathbb{R}^N}|2\pi\xi|^\alpha|\widehat{\psi}(\xi)|^2\,d\xi
    =\int_{\mathbb{R}^N}|(-\Delta)^{\alpha/4}\psi(x)|^2\, dx.\\
  \end{array}
$$

The extension operator is well defined for smooth functions through a Poisson kernel, whose explicit expression is given in~\cite{Caffarelli-Silvestre}.
 It  can also be defined in the space $\dot H^{\alpha/2}(\mathbb{R}^N)$, and in fact
\begin{equation}
  \label{eq:ext}
    \|\ext_\alpha(\psi)\|_{X^\alpha}=c_\alpha\|\psi\|_{\dot{H}^{\alpha/2}}\, , \quad \forall \, \psi\in \dot H^{\alpha/2}(\mathbb{R}^N),
\end{equation}
where $c_\alpha=\sqrt\kappa_\alpha$ , see Lemma~\ref{lemma:well-ext}.
On the other hand, for a function $\phi\in X^\alpha(\mathbb{R}^{N+1}_+)$, we will
denote its trace on $\mathbb{R}^N\times\{y=0\}$ as $\tr(\phi)$. This  trace operator is also well defined and it satisfies
\begin{equation}
  \label{eq:well-trace}
  \|\tr(\phi)\|_{\dot{H}^{\alpha/2}}\le
c_\alpha^{-1}\|\phi\|_{X^\alpha}.
\end{equation}

\subsection{A trace inequality}

In order to prove regularity of solutions to
problem~\eqref{prob:main.nonlocal} we will use a trace immersion. As
a first step we show that the corresponding result for the whole
space holds. Although most of the  results used in order to prove
Theorem~\ref{th:trace} below are known we have collected them for
the readers convenience.

First of all we prove inequality~\eqref{eq:well-trace}. The
Sobolev embedding yields then, that the trace also belongs to $L^{2^*_\alpha}(\mathbb{R}^N)$,
where $2^*_\alpha=\frac{2N}{N-\alpha}$. Even the best constant
associated to this inclusion is attained and can be characterized.

\begin{thm}\label{th:trace}
For every $z\in X^\alpha(\mathbb{R}^{N+1}_+)$ it holds
\begin{equation}
  \label{eq:trace}
  \left(\int_{\mathbb{R}^{N}} |v(x)|^{\frac{2N}{N-\alpha}}\,dx\right)^{\frac{N-\alpha}N}\le S(\alpha,N)
  \int_{\mathbb{R}^{N+1}_+} y^{1-\alpha}|\nabla z(x,y)|^2\,dxdy,
\end{equation}
where $v=\tr(z)$. The best constant takes the exact value
\begin{equation}
S(\alpha,N)=\frac{\Gamma(\frac\alpha2)\Gamma(\frac{N-\alpha}2)(\Gamma(N))^{\frac\alpha
N}} {2\pi^{\frac\alpha2}\Gamma(\frac{2-\alpha}2)
\Gamma(\frac{N+\alpha}2)(\Gamma(\frac N2))^{\frac\alpha
N}}
\label{constant-s}\end{equation} and it is achieved when $v$ takes
the form
\begin{equation}
v(x)=\big(|x-x_0|^2+\tau^2\big)^{-\frac{N+\alpha}2},
\label{func-min}
\end{equation} for some $x_0\in\mathbb{R}^N$, $\tau\in\mathbb{R}$, and $z=\ext_\alpha(v)$.
\end{thm}

The analogous results for the classical  Laplace operator can be found
in~\cite{Escobar,Lions85b}.

The proof of Theorem \ref{th:trace} follows from a series of lemmas, that we prove next.

\begin{lemma}
\label{lemma:well-ext}
    Let $v\in\dot H^{\alpha/2}(\mathbb{R}^N)$ and let $z=\ext_\beta(v)$ be its $\beta$-harmonic extension, $\beta\in(\alpha/2,2)$. Then
      $z\in X^{\alpha}(\mathbb{R}^{N+1}_+)$ and moreover there exists a
      positive universal constant $c(\alpha,\beta)$ such that
    \begin{equation}
      \|v\|_{\dot H^{\alpha/2}}=c(\alpha,\beta)\|z\|_{X^{\alpha}}.
    \label{eq:x2}\end{equation}
In particular if $\beta=\alpha$ we have
$c(\alpha,\alpha)=1/\sqrt{\kappa_\alpha}$.\end{lemma}
Inequality~\eqref{eq:trace} needs only the case $\beta=\alpha$,
which was included in the proof of  the local characterization of
$(-\Delta)^{\alpha/2}$ in~\cite{Caffarelli-Silvestre}. The
calculations performed in~\cite{Caffarelli-Silvestre} can be
extended to cover the  range $\alpha/2<\beta<2$ and in particular
includes the case $\beta=1$ proved in~\cite{Xiao}.

\demo Since $z=\ext_\beta(v)$, by definition $z$ solves
$\diver(y^{1-\beta}\nabla z)=0$,  which is equivalent to
$$
\Delta_x z+\frac{1-\beta}y\,\frac{\partial z}{\partial y}+\frac{\partial^2 z}{\partial
y^2}=0.
$$
Taking Fourier transform in $x\in\mathbb{R}^N$ for $y>0$ fixed, we have
$$
-4\pi^2|\xi|^2\hat z+\frac{1-\beta}y\,\frac{\partial \hat z}{\partial
y}+\frac{\partial^2
\hat z}{\partial y^2}=0.
$$
and $\hat z(\xi,0)=\hat v(\xi)$. Therefore $\hat z(\xi,y)=\hat
v(\xi)\phi_\beta(2\pi|\xi|y)$, where $\phi_\beta$ solves the problem
\begin{equation}
  -\phi+\frac{1-\beta}s
\phi'+\phi''=0,\qquad\phi(0)=1,\quad\lim\limits_{s\to\infty}\phi(s)=0.
\label{prob:phialpha}
\end{equation}
In fact, $\phi_\beta$ minimizes the functional
$$
H_\beta(\phi)=\int_0^\infty(|\phi(s)|^2+|\phi'(s)|^2)s^{1-\beta}\,ds.
$$
and it can be shown that it is a combination of Bessel functions,
see~\cite{Lebedev}. More precisely, $\phi_\beta$ satisfies the following asymptotic behaviour
\begin{equation}
  \begin{array}{ll}
\phi_\beta(s)\sim 1-c_1 s^{\beta},&\quad\mbox{for }\
s\to0,\\ [4mm] \phi_\beta(s)\sim c_2s^{\frac{\beta-1}2}e^{-s},&\quad\mbox{for
}\ s\to\infty, \end{array}
\label{bessel}\end{equation}
where
$$
c_1(\beta)=\frac{2^{1-\beta}\Gamma(1-\beta/2)}{\beta\Gamma(\beta/2)},\quad
c_2(\beta)=\frac{2^{\frac{1-\beta}2}\pi^{1/2}}{\Gamma(\beta/2)}.
$$

Now we observe that
$$
\begin{array}{rl}
\displaystyle\int_{\mathbb{R}^N} |\nabla z(x,y)|^2\,dx&=\displaystyle\int_{\mathbb{R}^N}
\left(|\nabla_x z(x,y)|^2+\Big|\frac{\partial z}{\partial
y}(x,y)\Big|^2\right)\,dx\\ [4mm]
&=\displaystyle\int_{\mathbb{R}^N}\left(4\pi^2|\xi|^2|\hat
z(\xi,y)|^2+\Big|\frac{\partial
\hat z}{\partial y}(\xi,y)\Big|^2\right)\,d\xi.
\end{array}$$
Then, multiplying by $y^{1-\alpha}$ and integrating in $y$,
$$
\begin{array}{l}
\displaystyle\int_0^\infty\int_{\mathbb{R}^N} y^{1-\alpha}|\nabla
z(x,y)|^2\,dxdy\\ [4mm] =\displaystyle
\int_0^\infty\int_{\mathbb{R}^N} 4\pi^2|\xi|^2|\hat
v(\xi)|^2\big(|\phi_\beta(2\pi|\xi|y)|^2+
|\phi_\beta'(2\pi|\xi|y)|^2\big)y^{1-\alpha}\,d\xi dy \\ [4mm]
=\displaystyle\int_0^\infty(|\phi_\beta(s)|^2+
|\phi_\beta'(s)|^2)s^{1-\alpha}\,ds\,\int_{\mathbb{R}^N}|2\pi\xi|^{\alpha}|\hat
v(\xi)|^2\,d\xi.
\end{array}$$
Using \eqref{bessel} we see that the integral
$\int_0^\infty(|\phi_\beta|^2+|\phi_\beta'|^2)s^{1-\alpha}\,ds$ is convergent
provided $\beta>\alpha/2$. This proves~\eqref{eq:x2} with
$c(\alpha,\beta)=(H_\alpha(\phi_\beta))^{-1/2}$.~\qed

\begin{obse}
If $\beta=1$ we have $\phi_1(s)=e^{-s}$, and
$H_\alpha(\phi_1)=2^{\alpha-1}\Gamma(2-\alpha)$, see~\cite{Xiao}. Moreover,
when $\beta=\alpha$, integrating by parts and using  the equation in~\eqref{prob:phialpha},
and~\eqref{bessel}, we obtain
\begin{equation}
\label{bessel2}
H_\alpha(\phi_\alpha)=\int_0^\infty
[\phi_\alpha^2(s)+(\phi_\alpha')^2(s)]s^{1-\alpha}\,ds=
-\lim_{s\to0}s^{1-\alpha}\phi_\alpha'(s)=\alpha c_1(\alpha)=\kappa_\alpha.
\end{equation}
\end{obse}

\begin{lemma}
Let $z\in X^\alpha(\mathbb{R}^{N+1}_+)$ and let
$w=\ext_\alpha(\tr(z))$ be its $\alpha$-harmonic associated function
(the extension of the trace). Then
  $$
    \|z\|_{X^\alpha}^2= \|w\|_{X^\alpha}^2+\|z-w\|_{X^\alpha}^2.
  $$
\label{lemma:minimum}\end{lemma}
  \demo Observe that, for $h=z-w$, we have
$$
\|z\|_{X^\alpha}^2=\int_{\mathbb{R}^{N+1}_+} y^{1-\alpha}(|\nabla
w|^2+|\nabla h|^2+2 \big\langle\nabla w,\nabla h\big\rangle).
$$
But, since $\tr(h)=0$, we have
$
\displaystyle\int_{\mathbb{R}^{N+1}_+} y^{1-\alpha}\big\langle\nabla
w,\nabla h\big\rangle\,dxdy =0
$.~\qed

\begin{lemma}
\label{lem:HL1}
 If $g\in L^{\frac{2N}{N+\alpha}}(\mathbb{R}^N)$, and $f\in \dot
H^{\alpha/2}(\mathbb{R}^N)$, then there exists a constant
$\ell(\alpha,N)>0$ such that
\begin{equation}
\left|\int f(x)\,g(x)\,dx\right|\le
\ell(\alpha, N)\|f\|_{\dot{H}^{\alpha/2}}\|g\|_{\frac{2N}{N+\alpha}}.
\label{fg-alpha}\end{equation}
Moreover, the equality in \eqref{fg-alpha} with the best constant holds when $f$ and $g$ take the
form~\eqref{func-min}.
\end{lemma}

The proof follows and standard argument that can be found, for instance in~\cite{cotsiolis, Xiao}.
\demo
By Par\c{c}eval's identity and Cauchy-Schwarz's inequality, we have
$$
\begin{aligned}
\displaystyle\left(\int_{\mathbb{R}^N} f(x)\,g(x)\,dx\right)^2&\displaystyle=
\left(\int_{\mathbb{R}^N} \widehat f(\xi)\,\widehat g(\xi)\,d\xi\right)^2 \\  &\displaystyle\le
\left(\int_{\mathbb{R}^N} |2\pi\xi|^\alpha\,|\widehat f(\xi)|^2\,d\xi\right)\,
\left(\int_{\mathbb{R}^N} |2\pi\xi|^{-\alpha}\,|\widehat g(\xi)|^2\,d\xi\right).
\end{aligned}$$
The second term can be written using~\cite{Lieb-Loss} as
$$
\int_{\mathbb{R}^N} |2\pi\xi|^{-\alpha}\,|\widehat g(\xi)|^2\,d\xi=b(\alpha,N)
\int_{\mathbb{R}^{2N}}\frac{g(x)g(x')}{|x-x'|^{N-\alpha}}\,dxdx',
$$
where
$$
b(\alpha,N)=\frac{\Gamma(\frac{N-\alpha}2)}{2^\alpha\pi^{N/2}\Gamma(\frac\alpha2)}\,,
$$
We now use the following Hardy-Littlewood-Sobolev inequality,
$$
\int_{\mathbb{R}^{2N}}\frac{g(x)g(x')}{|x-x'|^{N-\alpha}}\,dxdx'\le
d(\alpha,N)\|g\|^2_{\frac{2N}{N+\alpha}},
$$
see again~\cite{Lieb-Loss}, where
$$
d(\alpha,N)=\frac{\pi^{\frac{N-\alpha}2}\Gamma(\alpha/2)(\Gamma(N))^{\frac\alpha
N}}{\Gamma((N+\alpha)/2)(\Gamma(N/2))^{\frac\alpha N}},
$$
with equality if $g$ takes the form~\eqref{func-min}. From this we obtain the
desired estimate~\eqref{fg-alpha} with the constant $
    \ell(\alpha,N)=\sqrt{ b(\alpha,N)d(\alpha,N)}
$.

When applying Cauchy-Schwarz's inequality, we obtain an identity provided the
functions $|\xi|^{\alpha/2}\widehat f(\xi)$ and $|\xi|^{-\alpha/2}\widehat
g(\xi)$ are proportional. This means
$$
\widehat g(\xi)=c|\xi|^{\alpha}\widehat f(\xi)=c[(-\Delta)^{\alpha/2}
f]^{\wedge}\,(\xi).
$$
We end by observing that if $g$ takes the form~\eqref{func-min} and
$g=c(-\Delta)^{\alpha/2}f$ then $f$ also takes the form~\eqref{func-min}. In fact, the only
positive regular solutions to
$(-\Delta)^{\alpha/2}f=cf^{\frac{N+\alpha}{N-\alpha}}$ take the
form~\eqref{func-min}, see \cite{Chen-Li-Ou}.~\qed

\demodel{Theorem~{\rm\ref{th:trace}}}We apply Lemma~\ref{lem:HL1}
with $g=|f|^{\frac{N+\alpha}{N-\alpha}-1}f$, then use Lemma~\ref{lemma:well-ext} and conclude using
Lemma~\ref{lemma:minimum}. The best constant is
$S(\alpha,N)=\ell^2(\alpha,N)/\kappa_\alpha$.~\qed

\begin{obse}
If we let $\alpha$ tend to 2, when $N>2$, we recover the classical Sobolev
inequality for a function in $H^1(\mathbb{R}^N)$, with the same constant, see~\cite{talenti}. In order to pass to the limit in the right-hand side
of \eqref{eq:trace}, at least formally, we observe that $(2-\alpha)y^{1-\alpha}\,dy$ is a measure
on compact sets of $\mathbb{R}_+$ converging (in the weak-* sense) to a Dirac
delta. Hence
$$
\lim_{\alpha\to2^-}\int_0^1\left(\int_{\mathbb{R}^{N}}
|\nabla z(x,y)|^2\,dx\right)(2-\alpha)y^{1-\alpha}\,dy=\int_{\mathbb{R}^{N}}
|\nabla v(x)|^2\,dx.
$$
We then obtain
$$ \left(\int_{\mathbb{R}^{N}} |v(x)|^{\frac{2N}{N-2}}\,dx\right)^{\frac{N-2}N}\le S(N)
  \int_{\mathbb{R}^{N}} |\nabla v(x)|^2\,dx,
$$
with the best constant
$S(N)=\lim\limits_{\alpha\to2^-}\frac{S(\alpha,N)}{2-\alpha}=\frac{1}{\pi
N(N-2)}\Big(\frac{\Gamma(N)}{\Gamma(N/2)}\Big)^{\frac N2}$. It is achieved when
$v$ takes the form
\eqref{func-min} with $\alpha$ replaced by 2.
\end{obse}

\section{The fractional Laplacian in a bounded domain}
\label{sect:bounded}
\setcounter{equation}{0}

\subsection{Spectral decomposition}

To define the fractional Laplacian in a bounded domain we
follow~\cite{Cabre-Tan}, see also~\cite{stinga-torrea}. To this aim  we
consider the cylinder
$$
  \mathcal{ C}_\Omega=\{(x,y) : x\in\Omega,\  y\in
  \mathbb{R}_+\}\subset\mathbb{R}^{N+1}_+\,,
$$
and denote by $\partial_L\mathcal{ C}_\Omega$ its lateral boundary. We
also define the energy space
$$
X_0^\alpha(\mathcal{ C}_\Omega)=\{z\in L^{2}(\mathcal{ C}_\Omega)\,:\, z=0 \mbox{ on }
\partial_L\mathcal{ C}_\Omega,\  \displaystyle\int_{\mathcal{ C}_\Omega} y^{1-\alpha}|\nabla z(x,y)|^2\,dxdy<\infty\},
$$
with norm
$$
\|z\|_{X_0^\alpha}^2:=\int_{\mathcal{ C}_\Omega} y^{1-\alpha}|\nabla
z(x,y)|^2\,dxdy.
$$
We want to characterize the space $\Theta(\Omega)$, which is  the image of $X_0^\alpha(\mathcal{
C}_\Omega)$ under the trace operator,
$$
\Theta(\Omega)=\{u=\tr(w)\,:\, w\in X_0^\alpha(\mathcal{ C}_\Omega)\}.
$$
In particular we will show that the
fractional Laplacian in a bounded domain $\Omega$ is well defined for functions
in $\Theta(\Omega)$. To this aim, let us start considering the extension
operator and fractional Laplacian for smooth functions.

 \begin{defi}
Given a regular function $u$, we define its
$\alpha$-harmonic extension $w=\ext_\alpha(u)$ to the cylinder
$\mathcal{ C}_\Omega$ as the solution to the problem
\begin{equation}
  \left\{
    \begin{array}
      {rcl@{\qquad}l}
        -\diver(y^{1-\alpha}\nabla w)&=&0&\mbox{ in }\ \mathcal{ C}_\Omega,\\ [5pt]
        w&=&0&\mbox{ on }\  \partial_L\mathcal{ C}_\Omega,\\ [5pt]
        \tr(w)&=&u&\mbox{ on }\ \Omega.
    \end{array}
    \right.
\label{eq:local.prob.CS.bdd}\end{equation}
\end{defi}

As in the whole space, there is also a Poisson formula for the extension operator in a bounded domain, defined through the Laplace transform and the heat semigroup generator
$e^{t\Delta}$, see~\cite{stinga-torrea} for details.

\begin{defi}
The fractional operator $(-\Delta)^{\alpha/2}$ in $\Omega$, acting
on a regular function $u$, is defined by
\begin{equation}
    \label{normalder2}
(-\Delta)^{\alpha/2}u(x)=-\frac{1}{\kappa_\alpha}\lim\limits_{y\to0^+}y^{1-\alpha}\dfrac{\partial
w}{\partial y}(x,y),
\end{equation}
where $w=\ext_\alpha(u)$ and $\kappa_\alpha$ is given as in~\eqref{normalder}.
\end{defi}

It is classical that the powers of a positive operator in a bounded domain are
defined through the spectral decomposition using the powers of the eigenvalues
of the original operator. We show next that in this case this is coherent with
the Dirichlet-Neumann operator defined above.

\begin{lemma}
Let $(\varphi_j,\lambda_j)$ be the eigenfunctions and eigenvectors
of $-\Delta$ in $\Omega$ (with  Dirichlet  boundary data). Then $(\varphi_j,\lambda_j^{\alpha/2})$
are the eigenfunctions and eigenvectors of
$(-\Delta)^{\alpha/2}$ in $\Omega$ (with the same  boundary conditions).
Moreover,
 $\ext_\alpha(\varphi_j)=\varphi_j(x)\psi(\lambda_j^{1/2} y)$, where $\psi$
solves the problem
$$
 \left\{\begin{array}
      {rcl@{\qquad}l}
      \psi''+\dfrac{(1-\alpha)}s\psi'&=&\psi,&s>0,\\ [3mm]
      -\lim\limits_{s\to0^+}s^{1-\alpha}\psi'(s)&=&\kappa_\alpha,\\ [3mm]
      \psi(0)&=&1.
    \end{array}\right.
$$
\label{lem:spectral}\end{lemma}

The proof of this result, based on separating variables, is straightforward.
The function $\psi$ coincides with the solution $\phi_\alpha$ in
problem~\eqref{prob:phialpha}.

\begin{lemma}
\label{lemma:extension-formula}
  Let $u\in L^2(\Omega)$. If $u=\sum a_j
  \varphi_j,$ where $\sum a_j^2\lambda_j^{\alpha/2}<\infty$,
then  $\ext_\alpha(u)\in X_0^\alpha(\mathcal{C}_\Omega)$ and
$$
\ext_\alpha(u)(x,y)=\sum a_j\varphi_j(x)\psi(\lambda_j^{1/2} y).
  $$
\end{lemma}

\demo
The formula for the extension follows immediately from Lemma~\ref{lem:spectral}.

Put $w=\ext_\alpha(u)$. Using the orthogonality of the family
$\{\varphi_j\}$, together with $\int_\Omega\varphi_j^2=1$,
$\int_\Omega|\nabla\varphi_j|^2=\lambda_j$, and~\eqref{bessel2}, we
have
$$
\begin{array}{l}
\displaystyle\int_{\mathbb{R}^{N+1}_+}y^{1-\alpha}|\nabla
w(x,y)|^2\,dxdy\\ [4mm]
 = \displaystyle\int_0^\infty
y^{1-\alpha}\!\!\int_\Omega\!\!\Big( \sum
a_j^2|\nabla\varphi_j(x)|^2\psi(\lambda_j^{1/2}y)^2+
a_j^2\lambda_j\varphi_j(x)^2(\psi'(\lambda_j^{1/2}y))^2\Big)\,dxdy\\
[4mm] = \displaystyle\int_0^\infty y^{1-\alpha}\sum a_j^2\lambda_j\Big(
\psi(\lambda_j^{1/2}y)^2+ (\psi'(\lambda_j^{1/2}y))^2\Big)\,dy\\ [4mm] =
\displaystyle\sum a_j^2\lambda_j^{\alpha/2}\int_0^\infty s^{1-\alpha}\Big(
\psi(s)^2+ (\psi'(s))^2\Big)\,ds= \kappa_\alpha\displaystyle\sum
a_j^2\lambda_j^{\alpha/2}<\infty.
\end{array}
$$
\qed

As in Section~\ref{sect:whole-space}, we get that the extension operator minimizes the norm in $X_0^\alpha(\mathcal{ C}_\Omega)$.

\begin{lemma}
\label{lemma:min2}
Let $z\in X_0^\alpha(\mathcal{ C}_\Omega)$ and let
$w=\ext_\alpha(\tr(z))$. Then
  $$
    \|z\|_{X_0^\alpha(\mathcal{ C}_\Omega)}^2= \|w\|_{X_0^\alpha(\mathcal{ C}_\Omega)}^2+\|z-w\|_{X_0^\alpha(\mathcal{ C}_\Omega)}^2.
  $$
\end{lemma}

Lemmas~\ref{lemma:extension-formula} and~\ref{lemma:min2} imply that the space $\Theta(\Omega)$ coincides with the space
$$H^{\alpha/2}_0(\Omega)=\left\{u=\sum a_j
  \varphi_j\in L^2(\Omega)\  :\  \sum a_j^2\lambda_j^{\alpha/2}<\infty\right\},$$ equipped with the norm
$$
  \|u\|_{H^{\alpha/2}_0(\Omega)}=  \left(\sum a_j^2\lambda_j^{\alpha/2}\right)^{1/2}=\|(-\Delta)^{\alpha/4}u\|_{2}.
 $$
Observe that Lemma~\ref{lemma:extension-formula} gives
$$\|u\|_{H^{\alpha/2}_0(\Omega)}=\kappa_\alpha^{-1/2}\|\ext_\alpha(u)\|_{X^\alpha_0(\mathcal{C}_\Omega)}, $$
which is the same as in~\eqref{eq:ext} but for the spaces and norms involved.
As a direct consequence we get that the fractional Laplacian is well defined in
this space.

\begin{corol}
Let $u\in H^{\alpha/2}_0(\Omega)$, then $(-\Delta)^{\alpha/2}u=\sum
  a_j\lambda_j^{\alpha/2}\varphi_j$.
\end{corol}

Let now $v\in X_0^\alpha(\mathcal{ C}_\Omega)$. Its extension by zero
outside the cylinder $\mathcal{ C}_\Omega$ can be approximated by
functions with compact support in $\overline{\mathbb{R}^{N+1}_+}$.
Thus, the trace inequality~\eqref{eq:trace}, together with H\"older's
inequality, gives a trace inequality for bounded domains.

\begin{corol}\label{cor:trace-bdd}
For any $1\le r\le \frac{2N}{N-\alpha}$, and every $z\in
X_0^\alpha(\mathcal{ C}_\Omega)$, it holds
\begin{equation}
  \label{eq:trace-bdd}
  \left(\int_{\Omega} |v(x)|^{r}\,dx\right)^{2/r}\le C(r,\alpha,N,|\Omega|)
  \int_{\mathcal{ C}_\Omega} y^{1-\alpha}|\nabla z(x,y)|^2\,dxdy,
\end{equation}
where $v=\tr(z)$.
\end{corol}

\subsection{The linear problem}

We now use the extension problem~\eqref{eq:local.prob.CS.bdd} and
the expression~\eqref{normalder2} to reformulate the nonlocal
problems in a local way. Let $g$ be a regular function and consider
the following problems: the nonlocal problem
 \begin{equation}
  \label{prob:general}
\left\{
    \begin{array}{rcl@{\quad}l}
    (-\Delta)^{\alpha/2}u & = & g(x)&\mbox{ in }\  \Omega,\\
    u & = & 0&\mbox{ on }\  \partial \Omega,
    \end{array}
    \right.
    \end{equation}
and the corresponding local one
 \begin{equation}
  \left\{
    \begin{array}{rcl@{\quad}l}
        -\diver(y^{1-\alpha}\nabla w) & = & 0&\mbox{ in }\  \mathcal{ C}_\Omega,\\[5pt]
         w&=&0&\mbox{ on }\  \partial_L\mathcal{ C}_\Omega,\\ [5pt]
         -\dfrac1\kappa_\alpha\lim\limits_{y\to0^+}y^{1-\alpha}\dfrac{\partial w}{\partial
y}& = & g(x)&\mbox{ on }\  \Omega.
    \end{array}
    \right.
\label{prob:local.general}\end{equation}

We want to define the concept of solution to~\eqref{prob:general}, which is
done in terms of the solution to problem~\eqref{prob:local.general}.

\begin{defi}\label{def:weak.trace}
  We say that $w\in X^\alpha_0(\mathcal{ C}_\Omega)$ is an energy solution
to problem~\eqref{prob:local.general},  if for every
function $\varphi\in \mathcal{ C}_0^1(\mathcal{ C}_\Omega)$ it holds
\begin{equation}
    \int_{\mathcal{ C}_\Omega} y^{1-\alpha}\big\langle\nabla w(x,y),\nabla\varphi(x,y)\big\rangle\,dxdy=
    \int_\Omega \kappa_\alpha g(x)\varphi(x,0)\,dx.
    \label{energysol}
\end{equation}
\end{defi}
In fact more general test functions can be used in the above formula, whenever
the integrals make sense. A supersolution (subsolution) is a function that
verifies~\eqref{energysol} with equality replaced by $\ge$ ($\le$) for every
nonnegative test function.
\begin{defi}\label{def:weak}
  We say that $u\in H_0^{\alpha/2}(\Omega)$ is an energy solution to
problem~\eqref{prob:general} if $u=\tr(w)$ and $w$ is an energy
solution to problem~\eqref{prob:local.general}.
 \end{defi}

In order to deal with problem~\eqref{prob:local.general} we will assume,
without loss of generality, $\kappa_\alpha=1$, by changing the function $g$.

In~\cite{cabre-sire} this linear problem is also mentioned. There some results
are obtained using the theory of degenerate elliptic equations developed in~\cite{fabes-kenig-serapioni}, in particular a regularity result
for bounded solutions to this problem is obtained in~\cite{cabre-sire}. We
prove here that the solutions are in fact bounded if $g$ satisfies a minimal
integrability condition.

\begin{thm}\label{th:linear}
Let $w$ be a solution to problem~\eqref{prob:local.general}.  If $g\in
L^{r}(\Omega)$, $r>\frac{N}{\alpha}$, then $w\in L^\infty(\mathcal{ C}_\Omega)$.
\end{thm}
\demo The proof follows the same ideas as in~\cite[Theorem 8.15]{GT} and uses the
trace inequality \eqref{eq:trace-bdd}. Without loss of generality we may assume
$w\ge0$, and this simplifies notation. The general case is obtained in a
similar way.

We define for $\beta\geq 1$ and $K\geq k$ ($k$ to be chosen later) a
$\mathcal{ C}^1([k,\infty))$ function $H$, as follows:
$$
    H(z)=\left\{
    \begin{array}
      {l@{\qquad}l}
      z^\beta-k^\beta,&z\in[k,K],\\
      \mbox{linear},&z>K.
    \end{array}
    \right.
$$
Let us also define $v=w+k$ and choose as test function $\varphi$,
$$
  \varphi=G(v)=\int_k^v |H'(s)|^2\,ds,\qquad \nabla\varphi=|H´(v)|^2\nabla v.
$$
 Replacing this test function into the definition of energy solution we obtain on one hand:
 \begin{equation}
   \label{GT1}
     \begin{array}
   {rcl}
        \displaystyle\int_{\mathcal{ C}_\Omega} y^{1-\alpha}\big\langle\nabla w,\nabla\varphi\big\rangle\,dxdy &=&
    \displaystyle\int_{\mathcal{ C}_\Omega} y^{1-\alpha}\big|\nabla v|^2|H'(v)|^2\,dxdy
    \\[12pt]
     &=&\displaystyle \int_{\mathcal{ C}_\Omega} y^{1-\alpha}\big|\nabla H(v)|^2\,dxdy
     \\[12pt]
     &\geq& \displaystyle\bigg(\int_\Omega |H(v)|^{\frac{2N}{N-\alpha}} \,dx\bigg)^{\frac{N-\alpha}N}=
     \|H(v)\|^2_{\frac{2N}{N-\alpha}}\, ,
 \end{array}
  \end{equation}
 where the last inequality follows by \eqref{eq:trace-bdd}.

 On the other hand
 \begin{equation}
   \label{GT2}
  \begin{aligned}
   \int_\Omega g(x)\varphi(x,0)\,dx&=\int_\Omega g(x)G(v)\,dx\leq \int_\Omega g(x)vG'(v)\,dx\\
   &\leq \frac{1}{k}\int_\Omega g(x)v^2|H'(v)|^2\,dx= \frac{1}{k}\int_\Omega
   g(x)|vH'(v)|^2\,dx.
 \end{aligned}
 \end{equation}
Inequality~\eqref{GT1} together with~\eqref{GT2}, leads to
\begin{equation}\label{eq:ineq}
  \|H(v)\|_{\frac{2N}{N-\alpha}}\leq \left(\frac{1}{k}\|g\|_{r}\right)^{1/2}
  \|(vH'(v))^2\|^{1/2}_{\frac{r}{r-1}}=\|vH'(v)\|_{\frac{2r}{r-1}},
\end{equation}
by choosing $k=\|g\|_{r}$. Letting $K\to \infty$ in the definition
of $H$, the inequality~\eqref{eq:ineq} becomes
$$
  \|v\|_{\frac{2N\beta}{N-\alpha}}\leq \|v\|_{\frac{2 r\beta}{r-1}}.
$$
Hence for all $\beta\geq 1$ the inclusion $v\in L^{\frac{2
r\beta}{r-1}}(\Omega)$ implies the stronger inclusion $v\in
L^{\frac{2N\beta}{N-\alpha}}(\Omega)$, since
$\frac{2N\beta}{N-\alpha}>\frac{2 r\beta}{r-1}$ provided $r>\frac
N\alpha$. The result follows now, as in~\cite{GT}, by an iteration argument,
starting with $\beta=\frac{N(r-1)}{r(N-\alpha)}>1$ and $v\in
L^{\frac{2N}{N-\alpha}}(\Omega)$. This gives $v\in
L^\infty(\Omega)$, and then $w\in L^\infty(\mathcal{ C}_\Omega)$. In
fact we get the estimate
$$
\|w\|_\infty\le c(\|w\|_{X^\alpha}+\|g\|_r).
$$
\qed

\begin{corol}
\label{cor:GT}
  Let $w$ be a solution to problem~\eqref{prob:local.general}.
  If $g\in L^\infty(\Omega)$, then $w\in \mathcal{ C}^\gamma(\overline{\mathcal{ C}_\Omega})$ for some $\gamma\in(0,1)$.
\end{corol}

\demo
    It follows directly from Theorem~\ref{th:linear} and~\cite[Lemma 4.4]{cabre-sire}.
\qed

\section{The nonlinear nonlocal problem}\label{sect-f(u)}
\setcounter{equation}{0}

\subsection{The local realization}

We deal now with the core of the paper; i.e. the study of the nonlocal
problem~\eqref{prob:main.nonlocal}. We write that problem in local version in
the following way: a solution to problem~\eqref{prob:main.nonlocal} is a
function $u=\tr(w)$, the trace of $w$ on $\Omega\times\{y=0\}$, where $w$
solves the local problem

 \begin{equation}
  \left\{
    \begin{array}
      {rcl@{\qquad}l}
        -\diver(y^{1-\alpha}\nabla w)&=&0&\mbox{ in }\  \mathcal{ C}_\Omega,\\ [5pt]
        w&=&0&\mbox{ on } \partial_L\mathcal{ C}_\Omega,\\[5pt]
       \dfrac{\partial w}{\partial\nu^\alpha}&=&f(w)&\mbox{ in }\  \Omega,\\
        \end{array}
    \right.
\label{mainprob}\end{equation}
where
\begin{equation}\label{nualpha}
\frac{\partial
w}{\partial\nu^\alpha}(x)=-\lim\limits_{y\to0^+}y^{1-\alpha}\dfrac{\partial
w}{\partial y}(x,y).
\end{equation}

In order to simplify notation in what follows we will denote $w$ for the function defined in the cylinder $\mathcal{ C}_\Omega$ as well as for its trace $\tr(w)$ on $\Omega\times\{y=0\}$.

As we have said, we will focus on the particular nonlinearity
\begin{equation}\label{eq:cc}
  f(s)=f_\lambda(s)=\lambda s^q+s^p.
\end{equation}
However many auxiliary results will be proved for more general reactions $f$
satisfying the growth condition
\begin{equation}
  \label{eq:hyp.crecimiento.f}
  0\le f(s)\leq c(1+|s|^p), \qquad \mbox{for some\ } p>0.
\end{equation}

\begin{obse}
  \label{rem:kappa}
    In the definition \eqref{nualpha} we have neglected the constant $\kappa_\alpha$
appearing in~\eqref{normalder2} by a simple rescaling. Therefore,
the results on the coefficient $\lambda$ for the local
problem~\eqref{mainprob}--\eqref{eq:cc} in this section are
translated into problem~\eqref{prob:main.nonlocal} with $\lambda$
multiplied by $\kappa_\alpha^{p(q-1)-1}$.
\end{obse}

Following Definition~\ref{def:weak.trace}, we say that $w\in X_0^\alpha(\mathcal{
C}_\Omega)$ is an energy solution of~\eqref{mainprob} if the following identity
holds
$$ \int_{\mathcal{ C}_\Omega}y^{1-\alpha}\big\langle\nabla w,\nabla\varphi\big\rangle\,dxdy=
  \int_\Omega f(w)\varphi\,dx
$$
 for every regular test function $\varphi$. In the analogous way we define sub-
and supersolution.

We consider now the functional
$$ J(w)=\frac12\int_{\mathcal{ C}_\Omega} y^{1-\alpha}|\nabla w|^2\,dxdy-\int_{\Omega} F(w)\,
  dx,
$$
 where $F(s)=\int_0^s f(\tau)\,d\tau$. For simplicity of notation, we
define $f(s)=0$ for $s\le0$. Recall that the trace satisfies $w\in
L^{r}(\Omega)$, for every $1\le r\le \frac{2N}{N-\alpha}$ if
$N>\alpha$, $1<r\le\infty$ if $N\le\alpha$. In particular if $p\leq
\frac{N+\alpha}{N-\alpha}$, and $f$
verifies~\eqref{eq:hyp.crecimiento.f}  then $F(w)\in L^1(\Omega)$,
and the functional is well defined and bounded from below.

It is well known that critical points of $J$ are solutions to~\eqref{mainprob}
with a general reaction $f$. We consider also the minimization problem
$$
  I=\inf\Big\{\int_{{C}_\Omega} y^{1-\alpha}|\nabla w|^2\,dxdy\,:\,w\in X_0^\alpha({C}_\Omega),\,\int_{\Omega}F (w)\,
  dx=1\Big\} ,
$$
 for which, by classical variational techniques, one has that below the critical
exponent, the infimum $I$ is achieved.
\begin{prop}
  If $f$ satisfies~\eqref{eq:hyp.crecimiento.f} with $p < \frac{N+\alpha}{N-\alpha}$,
  then there exists a nonnegative function $w\in X_0^\alpha(\mathcal{ C}_\Omega)$
  for which $I$ is achieved.
\end{prop}

We now establish two preliminary results. The first one is a classical
procedure of sub- and supersolutions to obtain a solution. We omit its proof.

\begin{lemma}\label{lem:iteration}
Assume there exist a subsolution $w_1$ and a supersolution $w_2$ to
problem~\eqref{mainprob} verifying $w_1\le w_2$. Then there also exists a solution $w$
satisfying $w_1\le w\le w_2$ in~$\mathcal{ C}_\Omega$.
\end{lemma}

The second one is a comparison result for concave nonlinearities. The proof
follows the lines of the corresponding one for the Laplacian performed
in~\cite{Brezis-Kamin}.

\begin{lemma}\label{lem:comparison}
  Assume that $f(t)$ is a function verifying that $f(t)/t$ is decreasing for $t>0$. Consider $w_1,w_2\in  X^{\alpha}_0(\mathcal{C}_\Omega)$
  positive subsolution, supersolution respectively
  to problem~\eqref{mainprob}. Then $w_1\leq w_2$ in  $\overline{\mathcal{ C}_\Omega}$.
\end{lemma}

\demo By definition we have, for the nonnegative test functions
$\varphi_1$ and $\varphi_2$ to be chosen in an appropriate way,
$$
\begin{array}{l}
\displaystyle\int_{\mathcal{ C}_\Omega}y^{1-\alpha}\big\langle\nabla
w_1,\nabla\varphi_1\big\rangle\,dxdy\le
  \int_\Omega f(w_1)\varphi_1\,dx,  \\ [3mm]
\displaystyle\int_{\mathcal{ C}_\Omega}y^{1-\alpha}\big\langle\nabla
w_2,\nabla\varphi_2\big\rangle\,dxdy\ge
  \int_\Omega f(w_2)\varphi_2\,dx.
  \end{array}
$$
Now let $\theta(t)$ be a smooth nondecreasing function such that
$\theta(t)=0$ for $t\le0$, $\theta(t)=1$ for $t\ge1$, and set
$\theta_\varepsilon(t)=\theta(t/\varepsilon)$. If we put, in the
above inequalities
$$
\varphi_1=w_2\,\theta_\varepsilon(w_1-
w_2),\qquad\varphi_2=w_1\,\theta_\varepsilon(w_1-w_2),
$$
we get
$$ I_1\ge \displaystyle
  \int_\Omega w_1w_2\Big(\frac{f(w_2)}{w_2}-
  \frac{f(w_1)}{w_1}\Big)\theta_\varepsilon(w_1-w_2)\,dx,
   $$ where
$$
  I_1:= \displaystyle\int_{\mathcal{ C}_\Omega}y^{1-\alpha}\big\langle w_1\,\nabla
w_2-w_2\,\nabla w_1,\nabla(w_1-w_2)\big\rangle\,\theta_\varepsilon'(w_1-w_2)\,dxdy.
$$
Now we estimate $I_1$ as follows:
$$
\begin{array}{rl}
\displaystyle I_1&\displaystyle \le\int_{\mathcal{
C}_\Omega}y^{1-\alpha}\big\langle \nabla w_1,(w_1-w_2)
\nabla(w_1-w_2)\big\rangle\,\theta_\varepsilon'(w_1-w_2)\,dxdy\\
[3mm]
  &\displaystyle=\int_{\mathcal{
C}_\Omega}y^{1-\alpha}\big\langle \nabla
w_1,\nabla\gamma_\varepsilon(w_1-w_2)\big\rangle\,dxdy \end{array}$$
where $\gamma_\epsilon'(t)=t\theta_\varepsilon'(t)$. Therefore,
since $0\le \gamma_\varepsilon\le\varepsilon$, we have
$$
I_1\le\int_{\Omega}f(w_1)\gamma_\varepsilon(w_1-w_2)\,dx\le
c\varepsilon.
$$
We end as in~\cite{ABC}. Letting $\varepsilon$ tend to zero, we
obtain
$$
\int_{\Omega\cap\{w_1>w_2\}} w_1w_2\left(\frac{f(w_2)}{w_2}-
  \frac{f(w_1)}{w_1}\right)\,dx\le0,
$$
which together with the hypothesis on $f$ gives $w_1\le w_2$ in
$\Omega$. Comparison in $\mathcal{ C}_\Omega$ follows easily by the
maximum principle.~ \qed

Now we show that the solutions to problem
\eqref{mainprob}--\eqref{eq:hyp.crecimiento.f} are bounded and H\"older
continuous. Later on, in Section~\ref{sec-further}, we will obtain a uniform
$L^\infty$-estimate in the case where $f$ is given by \eqref{eq:cc} and the
convex power is subcritical.

\begin{prop}
  Let $f$ satisfy~\eqref{eq:hyp.crecimiento.f} with $p< \frac{N+\alpha}{N-\alpha}$, and let
  $w\in X_0^\alpha(\mathcal{ C}_\Omega)$ be an energy solution to problem~\eqref{mainprob}.
  Then $w\in L^\infty(\mathcal{ C}_\Omega)\cap \mathcal{ C}^{\gamma}(\Om)$ for some $0<\gamma<1$.
\end{prop}
\demo
The proof follows closely the technique of~\cite{Brezis-Kato}. As in the proof
of Theorem~\ref{th:linear}, we assume $w\ge0$. We consider, formally, the test
function $\varphi=w^{\beta-p}$, for some $\beta>p+1$. The justification of the
following calculations can be made substituting $\varphi$ by some approximated
truncature. We therefore proceed with the formal analysis. We get, using the
trace immersion, the inequality
$$
\left(\int_\Omega w^{\frac{(\beta-p+1)N}{N-\alpha}}\right)^{\frac{N-\alpha}N}\le
C(\beta,\alpha,N,\Omega)\int_\Omega w^{\beta}.
$$
This estimate allows to obtain the following iterative process
$$
\|w\|_{\beta_{j+1}}\le C\|w\|_{\beta_j}^{\frac{\beta_j}{\beta_j-p+1}},
$$
with $\beta_{j+1}=\frac{N}{N-\alpha}(\beta_j+1-p)$. To have
$\beta_{j+1}>\beta_j$ we need $\beta_j>\frac{(p-1)N}{\alpha}$. Since $w\in
L^{2^*_\alpha}(\Omega)$, starting with $\beta_0=\frac{2N}{N-\alpha}$, we get
the above restriction provided $p<\frac{N+\alpha}{N-\alpha}$. It is clear that
in a finite number of steps we get, for $g(x)=f(w(x,0))$, the regularity $g\in
L^r$ for some $r>\frac N\alpha$. As a consequence, we obtain the conclusion
applying Theorem~\ref{th:linear} and Corollary~\ref{cor:GT}.~\qed

\subsection{A nonexistence result}

\begin{thm}
Assume $f$ is a $C^1$ function with primitive $F$, and $w$ is an energy
solution to problem~\eqref{mainprob}. Then  the
following Pohozaev-type identity holds
$$
\frac12\int_{\partial_L\mathcal{ C}_\Omega}y^{1-\alpha}\big\langle (x,y),\nu\big\rangle|\nabla
w|^2\,d\sigma-N\int_\Omega F(w)\,dx+\frac{N-\alpha}2\int_\Omega wf(w)\,dx=0.
$$
\end{thm}
\demo
Just use the identity
$$
\begin{array}{l}
\displaystyle\big\langle (x,y),\nu\big\rangle y^{\alpha-1}\diver(y^{1-\alpha}\nabla w)
+\diver\Big[y^{1-\alpha}\Big(\big\langle(x,y),\nabla w\big\rangle-\frac12(x,y)|\nabla
w|^2\Big)\Big] \\ [4mm]
\displaystyle +\Big(\frac{N+2-\alpha}2-1\Big)|\nabla w|^2=0,
\end{array}
$$
where $\nu$ is the (exterior) normal vector  to $\partial\Omega$. It is
calculus matter to check this equality.
\qed

\begin{corol}
If $\Omega$ is starshaped and the nonlinearity $f$ satisfies the inequality
$((N-\alpha)sf(s)-2N F(s))\ge0 $, then  problem~\eqref{mainprob}  has no
solution. In particular, in the case $f(s)=s^p$ this means that there is no
solution for any $p\ge \frac{N+\alpha}{N-\alpha}$.
\label{supercritical}\end{corol}

The case $\alpha=1$ has been proved in~\cite{Cabre-Tan}. The
corresponding result for the Laplacian (problem
\eqref{prob:main.nonlocal} with $\alpha=2$) comes from
\cite{Pohozaev}.

\subsection{Proof of Theorems~\ref{thm:main}-\ref{thm:main2}}

We prove here Theorems~\ref{thm:main}-\ref{thm:main2} in terms of the solution of the local
version~\eqref{mainprob}. For the sake of readability we split the proof of Theorem~\ref{thm:main}  into
several lemmas.
From now on we will denote
$$
(P_\lambda)\equiv\left\{\begin{array} {rcl@{\qquad}l}
    -\diver(y^{1-\alpha}\nabla w)&=&0 &\mbox{ in }\  \mathcal{ C}_\Omega,\\ [5pt]
    w&=&0 &\mbox{ on }\  \partial_L\mathcal{ C}_\Omega,\\[5pt]
    \dfrac{\partial
w}{\partial\nu^\alpha} &=&\lambda w^q+w^p, \quad w>0 &\mbox{ in }\  \Omega,\\
\end{array}\right.
$$
and consider the associated energy functional
$$
J_\l (w) =
\displaystyle\frac12\int_{\mathcal{ C}_\Omega} y^{1-\alpha}|\nabla w|^2\,dxdy-
  \int_{\Omega} F_\l(w)\,dx,
$$
where
$$
F_\lambda(s)=\frac\lambda{q+1}s^{q+1}+\frac1{p+1}s^{p+1}.
$$

\begin{lemma}\label{lem:>}
Let $\Lan$ be defined by
$$
\Lan =\sup\{ \lan >0\ :\ \mbox{Problem
$({P}_{\lan })$ has solution}\} .
$$
Then $\Lan <\infty$.
\end{lemma}
\demo
Consider the eigenvalue problem associated to the first eigenvalue $\lambda_1$, and let
$\vfi_1 >0$ be an associated eigenfunction, i.e.,
$$
\left\{\begin{array}
{rcl@{\qquad}l} -\diver(y^{1-\alpha}\nabla \vfi_1)&=&0 &\mbox{ in }\  \mathcal{
C}_\Omega,\\ [5pt]
    \vfi_1&=&0 &\mbox{ on }\  \partial_L\mathcal{ C}_\Omega,\\[5pt]
   \dfrac{\partial
\vfi_1}{\partial\nu^\alpha}&=&\lambda_1 \vfi_1   &\mbox{ in }\  \Omega.
\end{array}\right.
$$
Then using $\vfi_1$ as a test function in $(P_\lan)$ we have that
\begin{equation}
  \label{eq:lem.>}
  \int_\Om (\lan w^q+w^p)\vfi_1\,dx=\lan_1\int_\Om w\vfi_1\, dx.
\end{equation}
There exist positive constants $c,\de$ such that $\lan
t^q+t^p>c\lan^\de t$, for any $t>0$. Since $u>0$ we obtain, using~\eqref{eq:lem.>}, that $c\lan^\de<\lan_1$
which implies $\Lambda<\infty$.~\qed

 This proves the third statement in Theorem~\ref{thm:main}.

\begin{lemma}\label{lem:<}
Problem $(P_\lan)$ has a positive solution for every $0<\lambda <\Lambda$.
Moreover, the family $\{w_\lambda\}$ of minimal solutions is increasing with
respect to~$\lambda$.
 \end{lemma}

\begin{obse}
  \label{rem:Lambda}
  Although this $\Lambda$ is not exactly the same as that of~Theorem~\ref{thm:main},
  see Remark~\ref{rem:kappa}, we have not changed the notation for the sake of simplicity.
\end{obse}
\demo
The associated functional to problem $(P_{\lambda})$  verifies, using
Corollary~\ref{cor:trace-bdd},
$$
\begin{array}{rcl}
  J_\l (w) & = & \displaystyle\frac12\int_{\mathcal{ C}_\Omega} y^{1-\alpha}|\nabla w|^2\,dxdy-
  \int_{\Omega} F_\l(w)\,
  dx \\ [3mm] & \ge &\displaystyle
\dfrac{1}{2}\int_{\mathcal{ C}_\Omega}y^{1-\alpha}|\nabla w|^2\,dxdy- \lambda
C_1\left(\int_{\mathcal{ C}_\Omega}y^{1-\alpha}|\nabla
w|^2\,dxdy\right)^{\frac{q+1}{2}}
\\ [3mm] & &\displaystyle -C_2\left(\int_{\mathcal{ C}_\Omega}
y^{1-\alpha}|\nabla w|^2\,dxdy\right)^{\frac{p+1}{2}},
\end{array}
$$
for some positive constants $C_1$ and $C_2$. Then for $\lambda$ small enough
there exist two solutions of problem  $(P_{\lambda})$, one given by
minimization and another one given by the Mountain-Pass Theorem,
\cite{AR}. The proof is standard, based on the geometry of the
function $g(t)=\frac12t^2-\lambda C_1t^{q+1}-C_2t^{p+1}$, see for
instance~\cite{GP} for more details.

We now show that there exists a solution for every $\lambda\in(0,\Lambda)$. By
definition of $\Lambda$, we know that there exists a solution corresponding to
any value of $\lambda$ close to $\Lambda$. Let us denote it by $\mu$, and let
$w_\mu$ be the associated solution. Now $w_\mu$ is a supersolution for all
problems $(P_\lambda)$ with $\lambda<\mu$. Take $v_\lambda$ the unique solution
to problem
\eqref{mainprob} with $f(s)=\lambda s^q$. Obviously $v_\lambda$ is a
subsolution to problem $(P_\lambda)$. By Lemma~\ref{lem:comparison}
$v_\lambda\le w_\mu$. Therefore by Lemma~\ref{lem:iteration} we conclude that
there is a solution for all $\lambda\in (0,\mu)$, and as a consequence, for the
whole open interval $(0,\Lambda)$. Moreover, this solution is the minimal one.
The monotonicity follows directly from the comparison lemma.
\qed

 This proves the first statement in Theorem~\ref{thm:main}.

\begin{lemma}\label{lem:=}
  Problem~$(P_\lan)$ has at least one solution if $\lambda=\Lambda$.
\end{lemma}
\demo
Let $\{\lambda_n\}$ be a sequence such that $\lambda_n\nearrow
\Lambda $. We denote by $w_n=w_{\lambda_n}$ the minimal solution to
problem $(P_{\lambda_n})$. As in~\cite{ABC}, we can prove that the linearized
equation at the minimal solution has nonnegative eigenvalues. Then  it follows,
as in~\cite{ABC} again, $J_{\lambda_n}(w_n)<0$. Since $J'(w_{\lambda})=0$, one
easily gets the bound $\|w_n\|_{X_0^\alpha(\mathcal{ C}_\Omega)}\le k$. Hence,
there exists a weakly convergent subsequence in $X_0^\alpha(\mathcal{ C}_\Omega)$
and as a consequence a weak solution of $(P_{\lambda})$ for $\lambda =\Lambda$.
\qed

 This proves the second statement in Theorem~\ref{thm:main} and finishes the proof of this theorem.

\demodel{Theorem~{\rm\ref{thm:main2}}}
In order to find a second solution we follow some arguments in~\cite{ABC}. It is essential to have that the
first solution is given as a local minimum of the associated functional,
$J_{\lambda}$. To prove this last assertion we modify some arguments developed in \cite{Alama}.

Let $\lambda_0\in (0,\Lambda)$ be fixed and consider
$\lambda_0<\bar{\lambda}_1<\Lambda$. Take $\phi_0=w_{\lambda_0}$, $\phi_1=w_{\bar\lambda_1}$ the two minimal
solutions to problem $(P_\l)$ with $\lambda=\lambda_0$ and
$\lambda=\bar{\lambda}_1$ respectively, then by comparison, $\phi_0<\phi_1$. We define
$$
M=\{w\in X^\alpha_0(\mathcal{ C}_\Omega):\,0\le w\le \phi_1\}.
$$
Notice that $M$ is a convex closed set of $X^\alpha_0(\mathcal{ C}_\Omega)$. Since $J_{\lambda_0}$ is
bounded from bellow in $M$ and it is semicontinuous on $M$, we get the
existence of $\underline{\omega}\in M$ such that $J_{\lambda_0}(\underline{\omega})=\inf_{w\in
M}J_{\lambda_0}(w)$.
Let $v_0$ be the unique positive solution to  problem
\begin{equation}\label{eq:concave}
\left\{
\begin{array}{rcl@{\qquad}l}
-\diver(y^{1-\alpha}\nabla v_0)&=&0 &\mbox{ in }\  \mathcal{ C}_\Omega,\\[5pt]
v_0 &=& 0 &\mbox{ on }\  \partial_L\mathcal{ C}_\Omega,\\[5pt]
\dfrac{\partial
v_0}{\partial\nu^\alpha} &=&  v_{0}^q &\mbox{ in }\  \Omega.
\end{array}\right.
\end{equation}
Since for $0<\epsi<<\lambda_0$, and
$J_{\lambda_0}(\epsi v_0)<0$, we have $\epsi v_0\in M$, then $\underline{\omega}\neq 0$.
Therefore $J_{\lambda_0}(\underline{\omega})<0$. By arguments similar to those in~\cite[Theorem 2.4]{Stru},
we obtain that $\underline{\omega}$ is a solution to problem~$(P_{\lambda_0})$. There are two possibilities:
\begin{itemize}
\item
If $\underline{\omega}\not \equiv w_{\lambda_0}$, then the result follows.

\item If  $\underline{\omega}\equiv w_{\lambda_0}$, we have
just to prove that  $\underline{\omega}$ is a local minimum of
$J_{\lambda_0}$. Then the conclusion  follows by using an argument close to the one in~\cite{ABC}, so we omit the complete details.
\end{itemize}
\noindent {\it We argue by contradiction}.

Suppose that $\underline{\omega}$
is not a local minimum of $J_{\lambda_0}$ in $X_0^\alpha (\mathcal{ C}_\Omega )$, then there exists a sequence
$\{v_n\}\subset X_0^\alpha (\mathcal{ C}_\Omega )$ such that $\|v_n-\underline{\omega}\|_{X_0^\alpha}\to 0$ and
$J_{\lambda_0}(v_n)<J_{\lambda_0}(\underline{\omega})$.

\noindent Let $w_n=(v_n-\phi_1)^{+}$ and
$z_n=\max\{0,\min\{v_n,\phi_1\}\}$.  It is clear that $z_n\in M$ and
$$
z_n(x,y)= \left\{
\begin{array}{l@{\quad}l}
 0 & \mbox{  if   }v_n(x,y)\le 0,\\[6pt]
  v_n(x,y) & \mbox{ if }0\le v_n(x,y)\le
\phi_1(x,y),\\[6pt]
  \phi_1(x,y) & \mbox{ if }\phi_1(x,y)\le v_n(x,y).
\end{array}
\right.
$$
We set
$$
\begin{array}{ll}
  T_n\equiv \{(x,y)\in \mathcal{ C}_\Omega :\,z_n(x,y)=v_n(x,y)\},\qquad &S_n\equiv
\text{supp}(w_n),\\[6pt]
\widetilde{T}_n=\overline{T_n}\cap\Om, \qquad &\widetilde{S}_n=S_n\cap\Om.
\end{array}
$$
Notice that $\text{supp}(v_n^+)=T_n\cup S_n$.
We claim that
\begin{equation}
  \label{claim}
 |\widetilde{S}_n|_{\Om}\to 0\quad\mbox{  as } n\to \infty,
\end{equation}
where $|A|_\Om\equiv
\int_{\Om}\chi_A (x)\, dx$.

\noindent By the definition of $F_\lambda$, we set
$F_{\lambda_0}(s)=\dfrac{\l_{0}}{q+1}s_+^{q+1}+\dfrac{1}{p+1}s_{+}^{p+1}$, for $s\in\mathbb{R}$,
and get
$$
\begin{aligned}
\dyle J_{\l_0}(v_n) =&
     \dyle  \frac12 \io \weig|\n v_n|^2\, dxdy -\int_\O F_{\lambda_0}(v_n)\,dx
\\
=& \dyle  \frac12 \int_{T_n} \weig|\n z_n|^2\, dx
dy-\int_{\widetilde{T}_n}F_{\lambda_0}(z_n)\,dx +\frac12 \int_{S_n}\weig  |\n v_n|^2\,dxdy\\
 &-\dyle\int_{\widetilde{S}_n}F_{\lambda_0}(v_n)\,dx+\frac12 \io \weig |\n v^-_n|^2\,dxdy\\
=&
\frac12 \dyle\int_{T_n} \weig|\n z_n|^2\, dxdy -\int_{\widetilde{T}_n}F_{\lambda_0}(z_n)\,dx\\
 & \dyle+\frac12 \int_{S_n} \weig |\n (w_n+\phi_1)|^2\,dxdy-\int_{\widetilde{S}_n}F_{\lambda_0}(w_n+\phi_1)\,dx\\
 &\dyle+\frac12 \io\weig |\n v^-_n|^2\, dxdy. \end{aligned}
$$
Since
$$
\io \weig|\n z_n |^2\, dxdy=\int_{T_n}\weig  |\n v_n |^2\,dxdy+\int_{S_n}\weig|\n\phi_1|^2\, dxdy
$$
and
$$
\int_\O F_{\lambda_0}(z_n)\,dx=\int_{\widetilde{T}_n}F_{\lambda_0}(v_n)\,dx+\int_{\widetilde{S}_n}
F_{\lambda_0}(\phi_1)\,dx,
$$
by using the fact that $\phi_1$ is a supersolution to $(P_\l)$
with $\l=\l_0$, we conclude that
\begin{eqnarray*}
J_{\l_0}(v_n)& = & J_{\l_0}(z_n)+\frac12 \int_{S_n}\weig  (|\n
(w_n+\phi_1)|^2-|\n \phi_1|^2)\,dxdy\\
&&-\int_{\widetilde{S}_n}(F_{\lambda_0}(w_n+\phi_1)-F_{\lambda_0}(\phi_1))\,dx+\frac12 \io\weig|\n
v^-_n|^2\,dxdy\\ & \ge &
J_{\l_0}(z_n)+\frac{1}{2}\|w_n\|_{X^\alpha_0}^2+\frac{1}{2}
\|v_n^-\|_{X^\alpha_0}^2\\
&&- \int_\O
\left\{F_{\lambda_0}(w_n+\phi_1)-F_{\lambda_0}(\phi_1)-(F_{\lambda_0})_{u}(\phi_1)w_n\right\}\,dx\\
&\ge&
J_{\l_0}(\underline{\omega})+\frac{1}{2}\|w_n\|_{X^\alpha_0}^2+\frac{1}{2}
\|v_n^-\|_{X^\alpha_0}^2\\
&&- \int_\O
\left\{F_{\lambda_0}(w_n+\phi_1)-F_{\lambda_0}(\phi_1)-(F_{\lambda_0})_{u}(\phi_1)w_n\right\}\,dx.
\end{eqnarray*}
On one hand, taking into account that $0<q+1<2$, one obtains that
$$0\le \frac{1}{q+1}(w_n+\phi_1)^{q+1}-\frac{1}{q+1}\phi_1^{q+1}-\phi_1^{q}w_n\le
\frac{q}{2}\frac{w_n^2}{\phi_1^{1-q}}.$$
The well known Picone's inequality (see \cite{Picone}) establish:
$$
|\n u|^2- \n \left(\dfrac{u^2}{v}\right)\cdot\n v\ge 0,
$$
for differentiable functions $v>0$, $u\ge 0$. In our case, by an approximation argument we get
$$
\l_0\int_\O \frac{w_n^2}{\phi_1^{1-q}}
dx \le \|w_n\|_{X_0^\alpha}^2.
$$
On the other hand, since $p+1>2$,
$$
\begin{array}{rcl}
0 & \le &\dyle \frac{1}{p+1}(w_n+\phi_1)^{p+1}-\frac{1}{p+1}\phi_1^{p+1}-\phi_1^{p}w_n \le
\frac{r}{2}w_n^2(w_n+\phi_1)^{p-1}\\ & & \\
 & \le &
C(p)(\phi_1^{p-1}w_n^2+w_n^{p+1}).
\end{array}
$$
Hence using that $p+1<2^*_\al$ and the claim~\eqref{claim}
$$
\int_\O\left\{ \frac{1}{p+1}(w_n+\phi_1)^{p+1}-\frac{1}{p+1}\phi_1^{p+1}-\phi_1^{p}w_n\right\}dx
\le o(1)\|w_n\|_{X_0^\alpha}^2. $$ As a consequence we obtain that
\begin{eqnarray*}
J_{\l_0}(v_n)&\ge&
J_{\l_0}(\underline{\omega})+\frac{1}{2}\|w_n\|_{X_0^\alpha}^2(1-q-o(1))+
\frac{1}{2}\|v_n^-\|_{X_0^\alpha}^2\\
&\equiv&
J_{\l_0}(\underline{\omega})+\frac{1}{2}\|w_n\|_{X_0^\alpha}^2(1-q-o(1))+o(1).
\end{eqnarray*}
Since $q<1$, there results that
$J_{\l_0}(\underline{\omega})>J_{\l_0}(v_n)\ge J_{\l_0}(\underline{\omega})$ for
$n>n_0$, a contradiction with the main hypothesis. Hence $\underline\omega $ is a minimum.

 To finish the proof we have to prove the claim~\eqref{claim}. For $\epsi>0$ small, and $\de>0$ ($\delta$ to be chosen later), we consider
$$
\begin{array}{rcl}
E_n & = & \{x\in
\Omega:\,v_n(x)\ge\phi_1(x)\quad \wedge\quad \phi_1(x)>\underline{\omega}(x)+\delta\},\\
 F_n & = & \{x\in
\Omega:\,v_n(x)\ge\phi_1(x)\quad\wedge\quad\phi_1(x)\le
\underline{\omega}(x)+\delta\}.
\end{array}
$$
Using the fact that
\begin{eqnarray*} 0 & = & |\{x\in \Omega:\,
\phi_1(x)< \underline{\omega}(x)\}|=\left|\bigcap_{j=1}^\infty \left\{x\in \Omega:\,
\phi_1(x)\le \underline{\omega}(x)+\frac{1}{j}\right\}\right|\\ & = &
\lim_{j\to \infty}\left|\left\{x\in \Omega:\, \phi_1(x)\le
\underline{\omega}(x)+\frac{1}{j}\right\}\right|,
\end{eqnarray*}
we get for $j_0$ large enough, that if $\delta<
\frac{1}{j_0}$ then  $$|\{x\in \Omega:\,
\phi_1(x)\le \underline{\omega}(x)+\delta\}|\le \frac{\varepsilon}{2}.$$ Hence
we conclude that $|F_n|_\Om\le \frac{\varepsilon}{2}$.

Since $\| v_n-\underline{\omega}\|_{X_0^\al}\to 0$ as $n\to \infty$, in particular by the trace embedding,
$\| v_n-\underline{\omega}\|_{L^2(\Omega)}\to 0$. We obtain that,
for $n\ge n_0$ large,
$$
\frac{\d^2\varepsilon}{2}\ge\io|v_n-\underline{\omega}|^2 dx\ge \int_{E_n}
|v_n-\underline{\omega}|^2 dx \ge \d^2|E_n|_{\Om}.
$$
Therefore $|E_n|_{\Om}\le \dfrac{\epsi}{2}$. Since $\widetilde{S}_n\subset F_n\cup E_n$
we conclude that $|\widetilde{S}_n|_\Om\le \epsi$ for $n\le n_0$. Hence
$|\widetilde{S}_n|_\Om\to 0$ as $n\to \infty$ and the claim follows.~\qed

\subsection{Proof of Theorems~\ref{thm:smallnorm}-\ref{th:universalbound} and further results}\label{sec-further}
We start with the uniform $L^{\infty}$-estimates for solutions to problem~\eqref{prob:main.nonlocal} in its local version given by $(P_\lan)$.
\begin{thm}\label{thm:uniform}
  Assume $\alpha\ge1$, $p<\frac{N+\alpha}{N-\alpha}$ and $N\geq 2$. Then there exists a constant $C=C(p,\Omega)>0$
  such that every solution to problem~$(P_\lan)$ satisfies
  $$
\|w\|_{\infty}\le C,
  $$
  for every $0\leq \lambda\leq\Lambda$.
\end{thm}

The proof is based in a rescaling method and the following two nonexistence results, proved in~\cite{dps}:

\begin{thm}
\label{thm:L1}
  Let $\alpha\ge1$. Then the problem in the half-space $\mathbb{R}^{N+1}_{+}$,
$$
  \left\{
    \begin{array}{rcl@{\quad}l}
        -\diver(y^{1-\alpha}\nabla w)&=&0&\mbox{ in }\mathbb{R}^{N+1}_{+}\\
       \dfrac{\partial w}{\partial\nu^\alpha}(x)&=&w^p(x,0)& \mbox{ on }\partial\mathbb{R}^{N+1}_{+}=\mathbb{R}^N
    \end{array}
    \right.
$$
  has no solution provided $p<\frac{N+\alpha}{N-\alpha}$, if $N\ge2$, or for every $p$ if $N=1$.
\end{thm}

\begin{thm}
    \label{thm:L2}
  The problem in the first quarter
  $$\mathbb{R}^{N+1}_{++}=
  \{z=(x',x_N,y)\,|\,x'\in\mathbb{R}^{N-1},\,x_N>0,\,y>0\},$$
$$
  \left\{
    \begin{array}{rcl@{\quad}l}
        -\diver(y^{1-\alpha}\nabla w)&=&0,&x_N>0,\,y>0,\\
       \dfrac{\partial w}{\partial\nu^\alpha}(x',x_N)&=&w^p(x',x_N,0),&\\
        w(x',0,y)&=&0,&
    \end{array}
    \right.
$$
  has no positive bounded solution provided $p<\frac{N+\alpha}{N-\alpha}$.
\end{thm}

\demodel{Theorem~{\rm \ref{thm:uniform}}}Assume by contradiction that there exists a sequence $\{w_n\}\subset X^{\al}_0(\mathcal{{C}}_\Om)$ of solutions to $(P_\lan)$ verifying that
$M_n=\|w_n\|_\infty\to\infty$, as $n\to \infty$. By the Maximum Principle, the maximum of $w_n$ is attained at a point $(x_n,0)$ where $x_n\in\Om$.
We define $\Om_n=\frac{1}{\mu_n}(\Om-x_n)$, with $\mu_n=M_n^{(1-p)/\alpha}$, i.e., we center at $x_n$ and dilate by $\frac{1}{\mu_n}\to\infty$ as
$n\to \infty$.

We consider the scaled functions
$$
v_n(x,y)=\dfrac{w_n(x_n+\mu_n x,\mu_n y)}{M_n}, \quad\mbox{for } x\in\Om_n,\: y\ge 0.
$$
It is clear that $\| v_n\|\le 1$, $v_n(0,0)=1$ and moreover
\begin{equation}
  \label{eq:vn}
  \left\{\begin{array}
{rcl@{\qquad}l}
-\diver(y^{1-\alpha}\nabla v_n)&=&0&\mbox{ in }\  \mathcal{ C}_{\Omega_n},\\[5pt]
    v_n&=&0&\mbox{ on }\  \partial_L\mathcal{ C}_{\Omega_n},\\[5pt]
\dfrac{\partial v_n}{\partial\nu^\alpha} &=&\lambda M_n^{q-p}v_n^q+v_n^p  &\mbox{ in }\  \Omega_n\times \{ 0\} .
\end{array}\right.
\end{equation}
By the Arzel\`a-Ascoli
Theorem, there exists a subsequence, denoted again by $v_n$, which converges to some function $v$ as $n\to\infty$. In order to see the problem satisfied by $v$ we pass to the limit in the weak formulation of~\eqref{eq:vn}. We observe that $\| v_n\|_\infty\le 1$ implies $\|v_n\|_{X_0^\alpha(\mathcal{C}_\Omega)}\leq C$.

We define $d_n=\dist (x_n,\p\Om)$, then there are two possibilities as $n\to\infty$ according the behaviour of the ratio $\frac{d_n}{\mu_n}$:
\begin{enumerate}
\item $\left\{\dfrac{d_n}{\mu_n}\right\}_n$ is not bounded.
\item $\left\{\dfrac{d_n}{\mu_n}\right\}_n$ remains bounded.
\end{enumerate}

In the first case, since $B_{d_n/\mu_n}(0)\subset \Om_n$, for another subsequence if necessary, it is clear that $\Omega_n$ tends to $\mathbb{R}^N$ and the limit function $v$ is a solution to
$$
\left\{\begin{array}
{rcl@{\qquad}l}
-\diver(y^{1-\alpha}\nabla v)&=&0\quad&\mbox{ in } \mathbb{R}_+^{N+1},\\[5pt]
\dfrac{\partial
v}{\partial\nu^\alpha}&=&v^p &\mbox{ on } \p\mathbb{R}_+^{N+1}.
\end{array}\right.
$$
Moreover, $v(0,0)=1$ and $v>0$ which is a contradiction with Theorem~\ref{thm:L1}.

In the second case, we may assume, again for another subsequence if necessary, that
$\frac{d_n}{\mu_n}\to s\ge 0$ as $n\to\infty.$ As a consequence, passing to the limit, the domains $\Om_n$ converge (up to a rotation) to some half-space $H_s=\{x\in\mathbb{R}^N :  x_N>-s\}$.    We obtain here that $v$ is a solution to
$$
\left\{\begin{array}
{rcl@{\qquad}l}
-\diver(y^{1-\alpha}\nabla v)&=&0\quad&\mbox{ in } H_s\times(0,\infty),\\[5pt]
\dfrac{\partial
v}{\partial\nu^\alpha}&=&v^p &\mbox{ on } H_s\times\{0\},
\end{array}\right.
$$
with $\|v\|_\infty=1$, $v(0,0)=1$. In the case $s=0$ this is a contradiction with the continuity of $v$. If $s>0$, the contradiction comes from
Theorem~\ref{thm:L2}.~\qed

We next prove Theorem~\ref{thm:smallnorm} in its local version.
\begin{thm}\label{lem:smallnorm}
 There exists at most one solution to problem $(P_\lan)$
with small norm.
\end{thm}
 We follow closely the arguments in \cite{ABC}, so we establish the following previous result:
\begin{lemma}
Let $z$ be the unique solution to problem \eqref{eq:concave}. There exists a constant $\be>0$ such that
\begin{equation}\label{eq:eigenpositive}
\|\phi\|^2_{X_0^\al (\mathcal{{C}}_\Om)}-q\int_{\Om}z^{q-1}\phi^2\, dx\ge \be\|\phi\|_{L^2(\Om)}^2,\quad \forall \phi\in X_0^\al(\mathcal{{C}}_\Om).
\end{equation}
\end{lemma}
\demo
We recall that $z$ can be obtained by minimization
$$
\min \left\{ \frac 12 \| \om\|^2_{X_0^\al (\mathcal{{C}}_\Om)}-\frac{1}{q+1}\| w\|^{q+1}_{L^{q+1}(\Om)}:\quad \om \in X_0^\al (\mathcal{{C}}_\Om)\right\}.
$$
As a consequence,
$$
 \| \phi \|^2_{X_0^\al (\mathcal{{C}}_\Om)}-q\int_\Om z^{q-1}\phi^2\, dx\ge 0,\quad \forall \phi\in X_0^\al(\mathcal{{C}}_\Om).
$$
This implies that the first eigenvalue $\rho_1$ of the linearized problem
$$
\left\{\begin{array}
{rcl@{\qquad}l}
-\diver(y^{1-\alpha}\nabla \phi)&=&0,\quad&\mbox{ in }\  \mathcal{ C}_\Omega,\\[5pt]
\phi&=&0,\quad&\mbox{ on }\  \partial_L\mathcal{ C}_\Omega,\\[5pt]
\dfrac{\partial
\phi}{\partial\nu^\alpha}-qz^{q-1}\phi&=&\rho\phi, &\mbox{ on }\  \Omega\times\{0\},
\end{array}\right.
$$
is nonnegative.

Assume first that $\rho_1=0$ and let $\varphi$ be a corresponding eigenfunction.
Taking into account that $z$ is the solution to~\eqref{eq:concave} we obtain that
$$q\int_\Om z^q\vfi\, dx=\int_\Om z^q\vfi\, dx$$ which is a contradiction.

Hence $\rho_1>0$, which proves \eqref{eq:eigenpositive}.
\qed

\demodel{Theorem {\rm \ref{lem:smallnorm}}}
Consider $A>0$ such that $pA^{p-1}<\be$, where $\be$ is given in
\eqref{eq:eigenpositive}.  Now we prove that problem $(P_\lan)$ has at most one
solution with $L^\infty$-norm less than $A$.

Assume by contradiction that  $(P_\lan)$ has a second solution $w=w_\lan+v$
verifying  $\| w\|_\infty<A$. Since $w_\lan$ is the minimal solution, it
follows that $v>0$ in $\Omega\times[0,\infty)$. We define now $\eta=\lan^{\frac{1}{1-q}}z$, where $z$ is the solution to
\eqref{eq:concave}. Then it verifies $-\div(\weig\n \eta)=0$, with boundary
condition $\lan \eta^q$. Moreover, $w_\lan$ is a supersolution to the
problem that $\eta$ verifies. Then by comparison, Lemma~\ref{lem:comparison},
applied with $f(t)=\lan t^q$, $v=\eta$ and $w=w_\lan$, we get
\begin{equation}\label{eq:ineq-concave}
w_\lan\ge \lan^{\frac{1}{1-q}}z \qquad \mbox{on } \Omega\times\{0\}.
\end{equation}
Since $w=w_\lan+v$ is solution to $(P_\lan)$ we have, on $\Omega\times\{0\}$,
$$
\dfrac{\partial (w_\lan+v)}{\partial\nu^\alpha}=\lan (w_\lan+v)^q+(w_\lan+v)^p\le
\lan w_\lan^q+\lan q w_\lan^{q-1}v+(w_\lan+v)^p,
$$
where the inequality is a consequence of the concavity, hence
$$
\dfrac{\partial v}{\partial\nu^\alpha}\le\lan q w_\lan^{q-1} v +(w_\lan +v)^p-w_\lan^p.
$$
Moreover, \eqref{eq:ineq-concave} implies $w_\lan^{q-1}\le \lan^{-1}z^{q-1}.$
From the previous two inequalities we get
$$
\dfrac{\partial v}{\partial\nu^\alpha}\le q z^{q-1} v +(w_\lan +v)^p-w_\lan^p.
$$
Using that $\|w_\lan +v\|_\infty\le A$, we obtain $(w_\lan+v)^p-w_\lan^p\le
pA^{p-1}v$. As a consequence,
$$
\dfrac{\partial v}{\partial\nu^\alpha}- q z^{q-1} v \le pA^{p-1}v.
$$
Taking $v$ as a test function and $\phi=v$ in \eqref{eq:eigenpositive} we arrive to
$$
\be\int_\Om v^2 \, dx\le p A^{p-1}\int_\Om v^2 \, dx.
$$
Since $pA^{p-1}<\be$ we conclude that $v\equiv 0$, which gives the desired
contradiction.~\qed
\begin{obse}
This proof also provides the asymptotic behavior of $w_\lan$ near $\lan= 0$,
namely $w_\lan \approx\lan^{ \frac{1}{1-q}}z$, where $z$ is the unique solution
to problem~\eqref{eq:concave}.
\end{obse}
\vskip .5cm
{\small
\noindent \textsc{Acknowledgments.}\\
 C.B. and A.dP.  are partially supported by
Spanish Project MTM2008-06326-C02-02; E.C. is partially supported by
Spanish Project MTM2009-10878.}

\end{document}